\def\date{04/06/2007}

\documentclass[envcountsame,envcountsect]{svmult}

\usepackage{makeidx}         
\usepackage{graphicx}        
\usepackage{multicol}        
\usepackage[bottom]{footmisc}
\usepackage{amsmath}
\usepackage{latexsym}
\usepackage{amssymb}

\makeindex             

\newcommand{\R}{{\mathbb R}}

\newcommand{\Z}{{\mathbb Z}}
\newcommand{\C}{{\mathbb C}}
\newcommand{\K}{{\mathbb K}}

\newcommand{\g}{{\mathfrak g}}
\newcommand{\gl}{{\mathfrak{gl}}}
\newcommand{\h}{{\mathfrak h}}
\newcommand{\q}{{\mathfrak q}}
\newcommand{\f}{{\mathfrak f}}
\newcommand{\e}{{\mathfrak e}}

\newcommand{\PE}{\mathop{\rm PE}\nolimits}

\newcommand{\Gl}{\mathop{\rm Gl}\nolimits}

\newcommand{\Hom}{\mbox{Hom}}
\newcommand{\Aut}{\mbox{Aut}}

\newcommand{\Der}{\mbox{\rm Der}}

\newcommand{\cA}{{\cal A}}

\newcommand{\XX}{{\cal X}}
\newcommand{\YY}{{\cal Y}}

\newcommand{\II}{{\cal I}}

\newcommand{\FF}{{\cal F}}
\newcommand{\GG}{{\cal G}}

\newcommand{\SSSS}{{\cal{S}}}
\newcommand{\Gras}{{\rm Gras}}
\newcommand{\Fon}{{\rm Fun}}
\newcommand{\Fun}{{\rm Fun}}

\newcommand{\Reg}{\mathop{\rm Reg}\nolimits}

\newcommand{\JJ}{{\cal{J}}}
\newcommand{\PP}{\Bbb{P}}

\newcommand{\PPP}{{\rm P}}
\newcommand{\SSS}{{\rm S}}

\newcommand{\ad}{{{\rm ad}}}

\newcommand{\msk}{\medskip}
\newcommand{\ssk}{\smallskip}
\newcommand{\nin}{\noindent}

\begin{document}

\title*{Jordan structures and non-associative geometry}

\author{Wolfgang Bertram} 
\institute{
Institut Elie Cartan, Universit\'{e} Nancy I, Facult\'{e}
des Sciences, B.P. 239, 54506 Vand\oe{}uvre-l\`{e}s-Nancy, Cedex, France
\texttt{bertram@iecn.u-nancy.fr}
}

\maketitle

{\small \noindent{\bf Abstract.\/}
We give an overview over constructions of geometries associated
to Jordan structures (algebras, triple systems and pairs), featuring 
analogs of these constructions with the Lie functor on the one hand and
with the approach of non-commutative geometry on the other hand.
}
\medskip

{\small \noindent{\bf Keywords:}
Jordan pair, Lie triple system,
graded Lie algebra, filtered Lie algebra, generalized
projective geometry, flag geometries, (non-) associative algebra and geometry
}

\smallskip
{\small \noindent{\bf AMS subject classification:}
17C37, 
17B70, 
53C15 
}

\section*{Introduction}

Let us compare two aspects of the vast mathematical topic
``links between
geometry and algebra'': on the one hand, the {\em Lie functor}
establishes a close relation between Lie groups (geometric side) and
Lie algebras (algebraic side); this is generalized by a correspondence
between {\em symmetric spaces} and {\em Lie triple systems} (see
\cite{Lo69}).
On the other hand, the philosophy of {\em Non-Commutative Geometry}
generalizes  the relation between usual,
geometric point-spaces $M$ (e.g., manifolds) and the commutative and
associative algebra $\Reg(M,\K)$ of ``regular'' (e.g., bounded,
smooth, algebraic,..., according to the context)
$\K$-valued functions
on $M$, where $\K=\R$ or $\C$, by replacing the algebra $\Reg(M,\K)$ by
more general, possibly non-commutative algebras $\cA$.
The interaction between these two aspects seems to be rather weak;
indeed, in some sense they are ``orthogonal'' to each other:
first, the classical setting of Lie's third theorem is finite-dimensional,
whereas the algebras $\cA$ of commutative or non-commutative geometry
are typically infinite dimensional;
second, and more importantly,
taking a commutative and associative algebra as input, 
we obtain as Lie bracket
$[x,y]=xy-yx=0$, and hence (if we forget the associative structure and
retain only the Lie bracket) we are left with constructing a Lie group
with zero Lie bracket. But this is not very interesting: 
it does not even capture the specific information
encoded in ``commutative geometries''.\footnote{Of course, this does not
exclude that {\em other} ways of associating Lie groups to function algebras
are interesting, for example by looking at their derivations and
 automorphisms;
but our point is precisely that such constructions 
need more input than
the trivial Lie bracket on the algebra of functions.}

\ssk
This remark suggests that, if one looks for a link between the two aspects
just mentioned, it should be interesting to ask for an analog
of the Lie functor for a class of algebras that contains faithfully
the class of commutative and associative algebras, but rather sacrifices
associativity than commutativity. 
At this point let us note that the algebras $\cA$ of non-commutative
geometry are of course always supposed to be {\it associative}, so that
the term ``associative geometry'' might be more appropriate than
``non-commutative geometry''. Indeed, behind the garden of associative
algebras starts the realm of general algebras, generically neither
associative nor commutative, and, for the time being, nobody has an idea
of what their ``geometric interpretation'' might be.
Fortunately, two offsprings of associative algebras grow not too far
away behind the garden walls:
{\it Lie algebras} right on the other side (the branches reach over the
wall so abundantly that some people even consider them as still belonging
to the garden), and {\it Jordan algebras}
a bit further. 
Let us just recall that the former are typically obtained by 
skew-symmetrizing an associative algebra, $[x,y]=xy-yx$, whereas the
latter are typically obtained by symmetrizing them, $x \bullet y =
\frac{1}{2}(xy+yx)$ (the factor $\frac{1}{2}$ being conventional, in
order to obtain the same powers $x^k$ as in the associative algebra).

\ssk
So let us look at {\it Jordan algebras} -- their advantage being that the
class of commutative associative algebras is faithfully embedded.
In this survey paper we will explain their geometric interpretation via
certain {\em generalized projective geometries},
emphasizing that this interpretation
 really combines both aspects mentioned above.
Of course, it is then important to 
include the case of {\em infinite-dimensional} algebras, and to
treat them in essentially the same way as finite-dimensional ones.
This is best done in a purely algebraic framework,
leaving aside all questions of topologizing our algebras and geometric spaces.
Of course, such questions form an interesting topic for the further 
development of the theory: in a very general setting (topological
algebraic structures over general
topological fields or rings) basic results are given in \cite{BeNe05},
and certainly many results from
the more specific setting of {\it Jordan operator algebras} (Banach-Jordan
structures over
$\K=\R$ or $\C$; see \cite{HOS84}) admit interesting geometric interpretations
in our framework. Beyond the Banach-setting, it should be interesting
to develop a theory of {\em locally convex topological Jordan
structures} and their geometries, taking up  the 
associative theory (cf.\ \cite{Bil04}).

\ssk
The general construction for Jordan algebras (and  for other members of
the family of Jordan algebraic structures, namely
{\it Jordan pairs} and {\it Jordan triple systems}, which in some sense
are easier to understand than Jordan algebras;
cf.\ Section 1) being explained in the
main text (Section 2 and Section 3), 
let us here just look at the special case of ``commutative
geometry'', i.e., the case of the commutative and associative
algebra $\cA = \Reg(M,\K)$ of
``regular $\K$-valued functions'' on a geometric space $M$
(here, $\K$ is a commutative field  with unit $1$; in the main text we
will allow also commutative rings with $1$).
Looking at $\cA$ as a Jordan algebra, we associate to it the space
$\XX := \Reg(M,\K \PP^1)$ of ``regular functions from $M$ into the projective
line $\K \PP^1$''. 
Of course, $\XX$ is no longer an algebra or a vector space; however, 
we can recover all these structures if we want: we call two functions $f,g
\in \XX$  ``transversal'', and write $f\top g$, 
 if $f,g \in \XX$ take different values at all points $p \in M$.
Now choose three functions $f_0,f_1,f_\infty \in \XX$ that are
mutually transversal. For each point $p \in M$, the value $f_\infty(p)$
singles out an ``affine gauge'' (by taking
$f_\infty(p)$ as point at infinity in the projective line $\K\PP^1$
and looking at the affine line $\K\PP^1 \setminus \{ f_\infty(p) \}$);
then $f_0(p)$ singles out a ``linear gauge'' (origin in the affine line)
and $f_1(p)$ a ``unit gauge'' (unit element in the vector line). 
Obviously, any such choice of $f_0$, $f_1$, $f_\infty$ leads to an 
identification of the set $(f_\infty)^\top$ of all functions that are
transversal to $f_\infty$ with the usual algebra $\cA$
of regular functions on $M$. In other words, $\cA$ and 
$(\XX;f_0,f_1,f_\infty)$ carry the same
information, and thus we see that
 the encoding via $\cA$ depends on various choices,
from which we are freed by looking instead at the ``geometric space'' $\XX$.
It thus becomes evident that $\cA$ really is a sort of ``tangent algebra''
of the geometric (``non-flat'') space $\XX$ at the point $f_0$,
in a similar way as the Lie algebra of a Lie group reflects its tangent
structure at the origin.  

\ssk
This way of looking at ordinary ``commutative geometry'' leads to
generalizations that are different from non-commutative geometry (in the
sense of A.\ Connes) but still have much in common with it.
For instance, both theories have concepts of ``states'', i.e., 
a notion that, in the commutative case, amounts to recover the point
space $M$ from the algebra of regular functions. 
Whereas in case of non-commutative geometry this concept relies
heavily on {\em positivity} (and hence on the ordered structure of the base
field $\R$), the corresponding concept for generalized projective geometries
is purely geometric and closely related to the notion of {\em inner ideals}
in Jordan theory (Section 4). 
There are many interesting open problems related to these items, some of them
mentioned in Section 4.


\ssk
Summing up, it seems  that the topic of geometrizing Jordan structures,
incorporating both ideas from classical Lie theory and  basic ideas from
associative geometry, is well suited for opening the problem of ``general
non-associative geometry'', i.e., the problem of finding geometries 
corresponding to more general non-associative algebraic 
structures. 

\bigskip
Acknowledgements. This paper is partly based on notes of a series of lectures 
given at the university of Metz in April 2007, and I would like to
thank Said Benayadi for inviting me at this occasion.


\ssk
Notation: Throughout, $\K$ is
 a commutative base ring with unit $1$ and such that $2$ and
$3$ are invertible in $\K$.

\section{Jordan pairs and graded Lie algebras}

\subsection{$\Z/(2)$-graded Lie algebras and Lie triple systems}

Let $(\Gamma,+)$ be an abelian group. A Lie algebra $\g$ over $\K$ is called
{\it $\Gamma$-graded} if it is of the form 
$$
\g = \bigoplus_{n \in \Gamma}
\g_n, \quad \mbox{with} \quad [\g_m,\g_n] \subset \g_{n+m}.
$$
Let us consider the example
$\Gamma =\Z/(2)$. Here, $\g=\g_{\overline 0} \oplus
\g_{\overline 1}$, with a subalgebra $\g_{\overline 0}$ and a subspace
$\g_{\overline 1}$ which is stable under $\g_{\overline 0}$ and such that
$[\g_{\overline 1}, \g_{\overline 1}] \subset   \g_{\overline 0}$.
It is then easily seen that
the linear map $\sigma:\g \to \g$ which is $1$ on
$\g_{\overline 0}$ and $-1$ on $\g_{\overline 1}$ is an involution of $\g$
(automorphism of order $2$). Conversely, every involution gives rise
to $\Z/(2)$-grading of $\g$, and hence such gradings
 correspond bijectively to {\it symmetric
Lie algebras} $(\g,\sigma)$, i.e., to Lie algebras $\g$
together with an automorphism $\sigma$ of order
$2$. Then $\g_{\overline 1}$ (the $-1$-eigenspace of $\sigma$) is 
stable under taking triple Lie brackets $[[x,y],z]$. 
This leads to the notion of
{\em Lie triple system}:

\begin{definition}
A {\em Lie triple system (LTS)} is a $\K$-module $\q$
together with a trilinear map 
$$
R: \q \times \q \times \q \to \q, \quad (X,Y,Z) \mapsto R(X,Y)Z =:[X,Y,Z]
$$
such that

\begin{description}
\item{(LT1)} $R(X,X)=0$
\item{(LT2)} $R(X,Y)Z+R(Y,Z)X+R(Z,Z)Y =0$ (the Jacobi identity)
\item{(LT3)} the endomorphism $D:=R(X,Y)$ is a {\em derivation} of the 
trilinear product, i.e.,
$$
D R(U,V,W)=R(DU,V,W)+R(U,DV,W)+R(U,V,DW).
$$
\end{description}
\end{definition}

\ssk
\nin If $\g$ is a $\Z/(2)$-graded Lie algebra, then $\g_{\overline 1}$ 
with $[X,Y,Z]:=[[X,Y],Z]$ becomes an LTS, and every LTS arises in this way
since one may reconstruct a Lie algebra from an LTS via the {\it
standard imbedding} (see [Lo69]):
let $\h \subset \Der(\q)$ the subalgebra of derivations of the LTS $\q$
spanned by all operators $R(X,Y)$, $X,Y \in \q$, and define a bracket
on $\g:=\q \oplus \h$ by
$$
[(X,D),(Y,E)]:=\big(DY-EX,[D,E] - R(X,Y) \big).
$$
One readily checks that $\g$ is a $\Z/(2)$-graded Lie algebra.
As a side remark, one may note that this construction (the {\em standard 
imbedding}) is in general not functorial -- see \cite{Sm05} for a
discussion of this topic.

\subsection{$3$-graded Lie algebras and Jordan pairs}

A Lie algebra is called {\it $2n+1$-graded} if it is $\Z$-graded, with
$\g_j=0$ if $j \notin \{ -n, -n+1,\ldots,n \}$.
The linear map $D:\g \to \g$ with $Dx=jx$ for $x \in \g_j$ is then a
derivation of $\g$; if it is an inner derivation, $D=\ad(E)$, the element
$E \in \g$ is called an {\em Euler operator}.

To any $\Z$-graded algebra one may associate a $\Z/(2)$-grading
$\g = \g_{even} \oplus \g_{odd}$, by putting together all homogeneous
parts with even, resp.\ odd index.
In the sequel we are mainly interested in $3$-graded and $5$-graded
Lie algebras; in both cases, $\g_{odd} = \g_1 \oplus \g_{-1}$.
We then let $V^\pm := \g_{\pm 1}$ and define trilinear maps by
$$
T^\pm:  V^\pm \times V^\mp \times V^\pm \to V^\pm, \quad
(x,y,z) \mapsto  [[x,y],z].
$$
The maps $T_\pm$ satisfy the identity

\begin{description}
\item{(LJP2)} 
$T^\pm(u,v,T^\pm(x,y,z))=
T^\pm (T^\pm(u,v,x),y,z) -
T^\pm(x,T^\mp(v,u,y),z) +$
\item{ }
${ }\quad \quad \quad \quad \quad \quad \quad \quad \quad
\quad \quad \quad \quad \quad \quad \quad
 T^\pm(x,y,T^\pm(u,v,z))$
\end{description}

\ssk
\nin Indeed, this is just another version of the identity
(LT3), reflecting the fact that $\ad[u,v]$ is a derivation of $\g$
respecting the grading since $[u,v] \in \g_0$. 
In the $3$-graded case, we have moreover, for all
$x,z \in V^\pm$, $y \in V^\mp$:

\begin{description}
\item{(LJP1)} 
$T^\pm(x,y,z)=T^\pm(z,y,x)$
\end{description}

\ssk
\nin This follows from the fact that
 $\g_1$ and $\g_{-1}$ are abelian, and so by the 
 Jacobi identity  
$[[x,y],z] - [[z,y],x] = [[z,x],y] = 0$.

\begin{definition}
A pair of $\K$-modules $(V^+,V^-)$ together with trilinear maps
$T^\pm: V^\pm \times V^\mp \times V^\pm \to V^\pm$ is called a 
{\em (linear) Jordan pair} if (LJP1) and (LJP2) hold.
\end{definition}

Every linear Jordan pair arises by the construction just described:
if $(V^+,V^-)$ is a linear Jordan pair, let
$\q:=V^+ \oplus V^-$ and define a trilinear map
$
\tilde T: \q^3 \to \q
$
via
$$
\tilde T\big( (x,x'),(y,y'),(z,z') \big) :=
\big(T^+(x,y',z),T^-(x',y,z') \big).
$$
Then the map $\tilde T$ satisfies the same identities as $(T^+,T^-)$,
with $T^+$ and $T^-$ both replaced by $\tilde T$ (in other words, $(\q,
\tilde T)$ is
a {\em Jordan triple system}, see below). Now we define a trilinear bracket
$\q^3 \to \q$ by
$$
\big[ (x,x'),(y,y'),(z,z') \big]:= \tilde T((x,x'),(y,y'),(z,z')) -
\tilde T((y,y'),(x,x'),(z,z')).
$$
By a direct calculation (cf.\ Lemma 1.4 below), one checks that $\q$
with this bracket is
a Lie triple system, and its standard imbedding
$\g = \q \oplus \h$ with $\g_{\pm 1} = V^\pm$ and
$\g_0=\h$, is a $3$-graded Lie algebra.
Without loss of generality we may assume that $\g$ contains an
Euler operator $E$: in fact, the endomorphism $E:\q \to \q$ which is
$1$ on $V^+$ and $-1$ on $V^-$, is a derivation of $\q$ commuting with
$\h$, and hence we may replace $\h$ by $\h + \K e$ in the construction
of the standard embedding.

For $5$-graded Lie algebras, the identity (LJP1) has te be replaced
by another, more complicated identity, which leads to the notion
of {\it Kantor pair}, see \cite{AF99}.
As for Jordan pairs, Kantor pairs give rise to Lie triple systems of the
form $\q =V^+ \oplus V^-$, where now the standard imbedding leads back
to a $5$-graded Lie algebra.

\subsection{Involutive $\Z$-graded Lie algebras}

An {\it involution of a $\Z$-graded Lie algebra} is an automorphism $\theta$
of order
$2$ such that $\theta(\g_j)=\g_{-j}$ for $j \in \Z$.
If $\g$ is $3$- or $5$-graded, we
 let $V:=\g_1$ and define 
$$
T:V \times V \times V \to V, \quad (X,Y,Z) \mapsto [[X,\theta Y],Z].
$$
Then $T$ satisfies the identity (LJT2) obtained from
(LJP2) by omitting the indices $\pm$, and 
if $\g$ is $3$-graded, we have moreover the analog of (LJP1).

\begin{definition}
 A $\K$-module
together $V$ with a trilinear map $T:V^3 \to V$
 satisfying the identities (JP1) and
(JP2) obtained from
(LJP1) and (LJP2) by omitting indices is called a
 {\it (linear) Jordan triple system} (JTS).
\end{definition}

 Every linear Jordan triple system
arises by the construction just described:
just let $V^+ := V^- :=V$ and $T^\pm := T$; then $(V^+,V^-)$ is a Jordan
pair carrying an involution (=isomorphism $\tau$ from $(V^+,V^-)$ onto the
Jordan pair $(V^-,V^+)$, namely $\tau(x,x')=(x',x)$).
Let $\g$ be the $3$-graded Lie algebra associated to this Jordan pair;
then the involution $\tau$ of the Jordan pair induces an involution
$\theta$ of $\g$.
By the way, these arguments show that
 Jordan triple systems are nothing but Jordan pairs with involution
(cf.\ [Lo75]). For any JTS $(V,T)$, the Jordan pair $(V^+,V^-)=(V,V)$ with
$T^\pm=T$ is called the {\it underlying Jordan pair}. 

\ssk
The Lie algebra $\g$ now carries {\em two} involutions, namely $\theta$ and the
involution $\sigma$ corresponding to the $\Z/(2)$-grading
$\g = \g_{even} \oplus \g_{odd}$.
Both involutions commute. In particular, $\sigma$ induces an involution of
the $\theta$-fixed subalgebra $\g^\theta$, with $-1$-eigenspace being
$\{ X + \theta(X) | \, X \in \g_1 \}$.
 The LTS-structure of this space
is described by the following lemma.

\begin{lemma}
{\rm (The Jordan-Lie functor.)}
If $T$ is a JTS on $V$, then
$$
[X,Y,Z ]= T(X,Y,Z) - T(Y,X,Z)
$$
defines a LTS on $V$. 
\end{lemma}

\begin{proof}
The proof of (LT1) is trivial, and (LT2) follows easily from the symmetry
of $T$ in the outer variables.
In order to prove (LT3), we define the endomorphism
$R(X,Y):=T(X,Y,\cdot) - T(Y,X,\cdot)$ of $V$.
Then $R(X,Y)$ is a derivation of the trilinear map $T$, as follows from
the second defining identity (JT2).
But every derivation of $T$ is also a derivation of $R$, since $R$ is simply
defined by skew-symmetrization of $T$ in the first two variables.
\end{proof}

The lemma
defines a functor from the category of JTS to
the category of LTS over $\K$, which we call the {\it Jordan-Lie functor}.
In general, 
it is neither injective nor surjective; all the more surprising is the fact
that, in the real, simple and finite-dimensional case, the Jordan-Lie
functor is not too far from setting up
a one-to-one correspondence  (classification results due to E.\ Neher, cf.\ 
tables given in \cite{Be00}). 

\subsection{The link with Jordan algebras }

Whereas all results from the preceding sections are naturally and easily
understood from the point of view of Lie theory, the results to presented
next are ``genuinly Jordan theoretic'' -- by this we mean that proofs by
direct calculation in $3$-graded Lie algebras are much less straightforward;
the reader may try to do so, in order to better appreciate the examples
to be given in the next section.

\ssk
In the following, $\g$ is a $3$-graded Lie algebra, with associated
Jordan pair $(V^+,V^-)=(\g_1,\g_{-1})$.
For $x \in V^-$, we define a $\K$-linear map
$$
Q^-(x):V^+ \to V^-, \quad y \mapsto Q^-(x)y:=\frac{1}{2} T^-(x,y,x) =
- \frac{1}{2} [x,[x,y]].
$$
In the same way $Q^+(y)$ for $y \in V^+$ is defined.
Then the following holds:

\msk
\nin {\bf (1) The Fundamental Formula.}
For all $x \in V^-$, $y \in V^+$,
$$
Q^-(Q^-(x)y) = Q^-(x)Q^+(y)Q^-(x).
$$

\msk
\nin {\bf (2) Meyberg's Theorem.} Fix $a \in V^-$. Then $V_a:=V^+$ with product
$$
x \bullet_a y:= \frac{1}{2} T^+(x,a,y) 
$$
is a Jordan algebra (called the {\em $a$-homotope algebra}).
 Recall (see, e.g., \cite{McC04}), that
a {\it (linear) Jordan algebra} is a $\K$-module $J$ with a bilinear and 
commutative product
$x \bullet y$ such that the identity

\begin{description}
\item{(J2)} $x \bullet (x^2 \bullet y) = x^2 \bullet (x \bullet y)$
\end{description}

\nin holds. 

\msk \nin {\bf (3) Invertible elements and unital Jordan algebras.}
An element $a \in V^-$ is called
{\em invertible in $(V^+,V^-)$} if the 
 the  operator 
$Q^-(a): V^+ \to V^-$
is invertible.
The element $a \in V^-$ is invertible in $(V^+,V^-)$
 if, and only if, the Jordan algebra
$V_a$ described in the preceding point admits a unit element; this unit
element is then
$a^\sharp:=(Q(a))^{-1}(a)$.
Moreover, the Jordan pair $(V^+,V^-)$ can then be recovered from
the Jordan algebra $(V,\bullet)=(V_a,\bullet_a)$ as follows:
let $V^+:=V^-:=V$ as $\K$-modules and let
$T:V^3 \to V$,
$$
T(x,y,z):= (x \bullet y) \bullet z - y \bullet (x \bullet z) + 
x \bullet (y \bullet z).
$$
Then $(V,T)$ is a Jordan triple system, and its underlying Jordan pair
$(V^+,V^-)$ is the one we started with.



\subsection{Some examples}


\msk \nin {\bf (1) Jordan pairs of rectangular matrices.}
Let $A$ and $B$ two $\K$-modules and $W:=A \oplus B$.
Let $I:W \to W$ be the linear map that is $1$ on $A$ and $-1$ on $B$. Then
the Lie algebra $\g:=\gl_\K(W)$ is $3$-graded, with
$\g_j = \{ X \in \g | \, [I,X]=j X \}$, $i=1,0,-1$.
In an obvious matrix notation, this corresponds to describing
$\g$ by $2 \times 2$-matrices, $\g_0$ as diagonal and
$\g_1$ as upper and $\g_{-1}$ as lower triangular matrices:
$$
\begin{pmatrix} \g_0 & V^+ \cr V^- & \g_0 \cr \end{pmatrix}
$$
Therefore
$$
(V^+,V^-)=(\Hom_\K(B,A),\Hom_\K(A,B))
$$
carries the structure of a linear Jordan pair.
The Jordan pair structure is given by
$$
T^\pm (X,Y,Z) = X Y Z + ZYX .
$$
In fact, this is proved by the following calculation
(using matrix notation for elements of $\g$):
$$
\Big[ \big[\begin{pmatrix} 0 & X \cr 0 & 0 \cr \end{pmatrix} ,
\begin{pmatrix} 0  & 0  \cr Y & 0 \cr \end{pmatrix} \big],
\begin{pmatrix} 0 & Z \cr 0 & 0 \cr \end{pmatrix} \Big] =
\big[\begin{pmatrix} XY & 0 \cr 0 & -YX  \cr \end{pmatrix} ,
\begin{pmatrix} 0  & Z  \cr 0 & 0 \cr \end{pmatrix} \big] =
\begin{pmatrix} 0 & XYZ + ZYX  \cr 0 & 0 \cr \end{pmatrix} .
$$
The proof of the
Fundamental Formula for this Jordan pair is easy since we may
calculate in the associative algebra $\Hom_\K(W,W)$, where
$Q(x)=\ell_x \circ r_x$ is simply composition of left and of right
translation by $x$. Similarly, the proof of Meyberg's Theorem 
becomes an easy exercise: for any $a \in \Hom_\K(W,W)$, the product
$(x,y) \mapsto xay$ is again an associative product, and by
symmetrizing and restricting to $V^+$, if $a \in V^-$, we get
the Jordan algebra $\bullet_a$.

\ssk
In general, the Jordan pair $(V^+,V^-)$ does not contain
invertible elements. For instance, if $\K$ is a field and
$A$ and $B$ are not isomorphic as vector spaces, then it is easily seen
that $Q^-(x)$ is never bijective. 

\ssk
Finally, it is clear that, if $A$ and $B$ are isomorphic as $\K$-modules,
and fixing such an isomorphism in order to
identify $A$ and $B$,
an element $a \in V^+=A=V^-$ is invertible if it is invertible 
in the usual sense in $\Hom(A,A)$, and then $\bullet_a$ has
$a^{-1}$ as unit element.
Note that, in particular, this happens for the Lie algebra
$\gl(2,\K)$, with $V^+ \cong V^- \cong \K$, where 
$x \bullet_a y = xay = yax$.


\msk \nin {\bf (2) Loop systems.}
As for any algebraic structure defined by identities, 
spaces of functions with values in such a structure, equipped 
with the ``pointwise product'', form again a structure of the given
type. In our case,
 if $M$ is a set and $(V^+,V^-)$ a Jordan pair, then
$(\Fun(M,V^+),\Fun(M,V^-))$  is a again a Jordan
pair, and similarly for Jordan triple systems and Jordan algebras.
Moreover, if $(V^+,V^-)$ corresponds to a $3$-graded Lie algebra $\g$,
then the function space corresponds to the
$3$-graded Lie algebra $\Fun(M,\g)$.
In particuler, ordinary functions
$(\Fun(M,\K),\Fun(M,\K))$ form a Jordan pair with pointwise product
$$
(T(f,g,h))(p) = 2 f(p) \, g(p) \, h(p).
$$
Under suitable assumptions, ``regular'' functions (continuous, smooth,
algebraic,...) will form subpairs.

\msk \nin {\bf (3) The classical examples.} Besides rectangular matrices,
these include:

\begin{description}
\item{$\bullet$} full asociative algebras $A$. Here $(V^+,V^-)=(A,A)$, and
$\g \subset \gl(2,A)$ is the subalgebra generated by the strictly
upper triangular
and strictly  lower triangular matrices;
\item{$\bullet$} Hermitian elements of an involutive associative algebra
$(A,*)$. Here,
 $(V^+,V^-)= (A^*,A^*)$, $\g$ is the symplectic algebra of $(A,*)$
(cf.\ \cite{BeNe04}, Section 8.2). Note that
 symmetric and Hermitian matrices are
a special case;
\item{$\bullet$} skew-Hermitian elements of an involutive associative algebra
$(A,*)$: similar as above, replacing $*$ by its negative. As above,
Jordan pairs of skew-symmetric or skew-Hermitian matrices are a special case;
\item{$\bullet$} conformal geometries (or ``spin factors''). Here,
$V^+ = V^- =V$ is a $\K$-module with a non-degenerate
symmetric bilinear form
 $(\cdot | \cdot):V \times V \to \K$ 
 and trilinear map 
$$
T(x,y,z)= (x | z)y - (x | y)z - (z| y)x \, .
$$
Then $\g$ is the orthogonal algebra of the quadratic form on
 $\K \oplus V \oplus \K$ given by 
$$
\beta((r,v,s),(r,v,s))= rs + (v|v) \, .
$$ 
\end{description}

\nin 
There are also exceptional Jordan systems, namely octonionic
$1 \times 2$-matrices, resp.\ the Hermitian octonionic 
$3 \times 3$-matrices, corresponding to $3$-gradings of the exceptional
Lie algebras of type $E_6$, resp.\ $E_7$.

\section{The generalized projective geometry of a Jordan pair}

\subsection{The construction}

The following construction of a pair $(\XX^+,\XX^-)$ of homogeneous spaces
associated to a Jordan pair $(V^+,V^-)$ is due to J.R.\ Faulkner 
and O.\ Loos (see \cite{Lo95}): 
starting with a (linear) Jordan pair $(V^+,V^-)$, let $\g$ be the
$3$-graded Lie algebra $\g$ constructed  in Section 1.2; as explained there,
we may assume that
 $\g$ contains an {\em Euler
operator} $E$, i.e., an element such that $[E,X]=iX$ for $X \in \g_i$.

\ssk
Next we define two abelian groups $U^\pm = \exp(\g_{\pm 1})$ by observing
that, for $X \in \g_{\pm 1}$, the operator $\ad(X):\g \to \g$
is $3$-step nilpotent, and that
$$
\exp(X):=e^{\ad(X)} = {\rm id} + \ad(X) + \frac{1}{2} \ad(X)^2
$$
is an automorphism of $\g$. Since $\g_{\pm 1}$ are abelian, $U^\pm$
are abelian subgroups of $\Aut(\g)$.
The maps $\exp_\pm:\g_\pm \to \Aut(\g)$ are injective
(here we use our Euler operator: $\exp_\pm(X)={\rm id}$ implies
$\exp_\pm(X)E=E \pm X=E$, whence $X=0$). Thus $U^\pm$ is
isomorphic to $(V^\pm,+)$.
The {\em elementary projective group} associated to the Jordan pair
$(V^+,V^-)$ is the subgroup
$$
G:=\PE(V^+,V^-):= \langle U^+,U^- \rangle
$$
of $\Aut(\g)$ generated by $U^+$ and $U^-$.
Let $H$ be the subgroup of $G$ stabilizing the grading (i.e., commuting with
$\ad(E)$); then $U^+$ and $U^-$ are normalized by $H$, and  the
groups $P^\pm$ generated by $H$ and $U^\pm$ are semidirect products:
$$
P^\pm := \langle H, U^\pm \rangle \cong H \rtimes V^\pm .
$$
Finally, we define two homogeneous spaces
$\XX^\pm := G / P^\mp$ with base points $o^\pm = e / P^\mp$.

\ssk
The pair of ``geometric spaces'' $(\XX^+,\XX^-)$ is called the
{\it generalized projective geometry associated to $(V^+,V^-)$}.
It is indeed the geometric object associated to the Jordan pair
$(V^+,V^-)$, in a similar way as a Lie group is the geometric object associated
to a Lie algebra. Let us explain this briefly.

\subsection{Generalized projective geometries}

The spaces $\XX^\pm$ being defined as above, 
the direct product $\XX^+ \times \XX^-$ is equipped with a {\em transversality
relation}: call $(x,\alpha) \in \XX^+ \times \XX^-$ {\em transversal}, and
write $x \top \alpha$, if they are conjugate under $G$ to the base point
$(o^+,o^-)$. It is easy to describe the set $\alpha^\top := \{ x \in \XX^+
| \, x \top \alpha \}$ : if $\alpha =o^-$ is the base point (whose 
stabilizer is $P^+$), then
this is the $P^+$-orbit $P^+.o^+$ which is isomorphic to $V^+ \cong U^+.o^+$.
Note here that the set $\alpha^\top$ carries a canonical
structure of an affine space over $\K$, since the vector group $U^+$ acts
on it simply transitively.
By homogenity, the same statements hold for any $\alpha \in \XX^-$.
This observation leads to the following definition:

\begin{definition}
A {\em pair geometry} is a pair of sets $\XX=(\XX^+,\XX^-)$ together
with a binary relation $\top \subset \XX^+ \times \XX^-$
(called {\em transversality}, and we write $x \top \alpha$ for
$(x,\alpha) \in \top$) such that
$\XX^\pm$ is covered by subsets of the form
$$
\alpha^\top = \{ x \in \XX^\pm | \, x \top \alpha \}, \quad \alpha \in
\XX^\mp.
$$
A {\em linear pair geometry (over $\K$)} is a
 pair geometry $(\XX^+,\XX^-, \top)$
such that, for any transversal pair $x \top \alpha$, the set
$\alpha^\top$ is equipped with a structure of linear space over $\K$
(i.e., a $\K$-module), with origin $x$, and the same property holds by
exchanging the r\^oles of $\XX^+$ and $\XX^-$.

An {\em affine pair geometry} is a linear pair geometry such that
the underlying
affine structure of the linear space $(\alpha^\top,x)$ does not depend on
$x$, for all transversal pairs $(x,\alpha)$, and dually.
In other words, for all $\alpha \in \XX^-$, the set $\alpha^\top$ carries
the structure of an affine space over $\K$, and dually.
\end{definition}

The discussion above shows that, to any Jordan pair $(V^+,V^-)$ over $\K$,
we can associate an affine pair geometry $(\XX^+,\XX^-,\top)$.
The link with triple products is given by introducting the
{\em structure maps of a linear pair geometry}:
if $x \top \alpha$, $y \top \alpha$, $z \top \alpha$ and $r \in \K$, then let
$r_{x,\alpha}(y):= r y$ denote the product $r \cdot y$, 
and $y+_{x,\alpha} z$ the sum of 
$y$ and $z$ in the $\K$-module
$\alpha^\top$ with zero vector $x$. In other words, we define maps of three
(resp.\ four) arguments by
\begin{eqnarray*}
\PPP_r :&
 (\XX^+ \times \XX^- \times \XX^+)^\top  \to \XX^+, & 
(x,\alpha,y)   \mapsto \PPP_r(x,\alpha,y) :=r_{x,\alpha}(y),
\cr
\SSS: &(\XX^+ \times \XX^- \times \XX^+ \times \XX^+)^\top  \to 
\XX^+,  &
(x,\alpha,y,z)   \mapsto \SSS(x,\alpha,y,z) := y +_{x,\alpha} z,
\end{eqnarray*}
where the domain of $\PPP_r$ is the ``space of generic triples'',
$$
 (\XX^+ \times \XX^- \times \XX^+)^\top  = 
\{ (x,\alpha,y) \in \XX^+ \times \XX^- \times \XX^+ \vert \,
x \top \alpha, \, y \top \alpha \},
$$
and the domain of $\SSS$ is the similarly defined space of generic quadruples.
Of course, we should rather write $\PPP_r^+$ instead of $\PPP_r$ and
$\SSS^+$ instead of $\SSS$; dually, the maps
 $\PPP_r^-$ and $\SSS^-$ are then also defined.
The structure maps encode all the information of a linear pair geometry:
by fixing the pair $(x,\alpha)$, the structure maps describe the linear
structure of $(\alpha^\top,x)$, resp.\ of $(x^\top,\alpha)$.
In this way, linear pair geometries can be regarded as algebraic objects
whose structure is defined by (one or several) ``multiplication maps'', 
just like groups, rings, modules, symmetric spaces...

\ssk
This point of view naturally leads to the question whether there are
more specific``identities'' satisfied by the structure maps.
If $(\XX^+,\XX^-;\top)$ is the geometry associated to a Jordan pair,
then this is indeed the case: there are two such identities,
 denoted by (PG1) and (PG2) in \cite{Be02}; the first identity (PG1)
can be seen as an ``integrated version'' of the defining identity
(LJP2) of a Jordan pair; it implies the existence of a ``big''
automorphism group, whereas (PG2) is rather a global version of
the Fundamental Formula and implies the existence of ``structural maps''
exchanging the two partners $\XX^+$ and $\XX^-$
 (we refer the reader to \cite{Be02} for details).
The important point about these identities is that we can also go the
other way round: by ``deriving'' them in a suitable way, one can show
that any affine pair geometry satisfying (PG1) and (PG2) in all scalar
extensions gives rise to a Jordan pair $(V^+,V^-)$ as ``tangent geometry''
with respect to a fixed base point $(o^+,o^-)$.
This construction clearly parallels the correspondence between Lie groups
and Lie algebras (and even more closely the one between symmetric spaces
and Lie triple systems, \cite{Lo69}; see also \cite{Be06}); 
however, in contrast
to Lie theory, the constructions are  algebraic in nature and therefore
work much more generally in arbitrary dimension and over general base
fields and -rings. They define
 a bijection (even an equivalence of categories)
between Jordan pairs and connected generalized projective geometries
with base point (and, in the same way, between Jordan triple systems,
resp.\ unital Jordan algebras, and {\em generalized polar (resp.\
null) geometries}, i.e., generalized projective geometries with the
additional structure of a certain kind of involution, cf.\ \cite{Be02},
\cite{Be03}) . 
Thus, in principle, all Jordan theoretic notions can be translated
into geometric ones; in general, it is not at all obvious what the 
correct translation should be, but once it has been found, it often
sheds new light onto the algebraic notion. In Chapter 4 we illustrate
this at the example of the notion of
 {\it inner ideals} and their {\it complementation
relation}.

\subsection{The geometric Jordan-Lie functor}

{\em Symmetric spaces} are well-known examples of ``non-associative
geometries'' (cf.\ \cite{Lo69}): they are manifolds $M$ together with point reflections
$\sigma_x$ attached to every point $x \in M$ such that the {\em
multiplication map} $\mu:M \times M \to M$, $(x,y) \mapsto \sigma_x(y)$
satisfies the properties 

\begin{description}
\item{(M1)} $\mu(x,x)=x$,
\item{(M2)} $\mu(x,\mu(x,y)) = y$,
\item{(M3)} $\mu(x,\mu(y,z))=\mu(\mu(x,y),\mu(x,z))$,
\item{(M4)} the tangent map $T_x(\sigma_x)$ is the negative of the identity
of the tangent space $T_x M$.
\end{description}

\nin 
Recall from Lemma 1.4 the (algebraic) Jordan-Lie functor, associating a 
Lie triple system to each Jordan triple system. 
The  geometric analog of this functor associates to each
generalized polar geometry $(\XX^+,\XX^-;p)$ the  {\em symmetric space}
$M^{(p)}$ of its non-isotropic points. 
The geometric construction is very simple: 
essentially, the multiplication map $\mu$ is obtained from the structure
maps $\PPP_r$ by taking the scalar $r=-1$. More precisely,
given a polarity
(i.e., 
a pair of bijections $(p^+,p^-):(\XX^+,\XX^-) \to (\XX^-,\XX^+)$
 that are inverses of each other, compatible with
the structure maps and such that there exists a {\em non-isotropic point}
$x \in \XX^+$, i.e., $x \top p^+(x)$), for each non-isotropic point $x$
we define the point-reflection by 
$$
\sigma_x(y) := \PPP_{-1}(x,p^+(x),y).
$$
The above mentioned identity (PG1) then implies in a straightforward way
that the set $M^{(p)}$ of non-isotropic points of $p$ becomes a symmetric
space (cf.\ \cite{Be02} for the purely algebraic setting and \cite{BeNe05}
for the smooth setting). Moreover, the Lie triple system of the symmetric
space $M^{(p)}$ is precisely the one defined by Lemma 1.4. 
The well-known construction of finite and infinite dimensional 
{\em bounded symmetric domains} (cf.\ \cite{Up85}) is a special case of
the geometric Jordan-Lie functor, but, as mentioned after Lemma 1.4, we
obtain in fact many other symmetric spaces in this way (essentially, all
classical finite-dimensional real ones and about half of the exceptional
ones, cf.\ tables in \cite{Be00}). 

\subsection{Examples revisited}

\nin {\bf (1) Rectangular matrices correspond to
Grassmannians.}
Let $W = A \oplus B$ be as in Example (1) of Section 1.5 and let
$\Gras_A^B(W)$
be the set of all submodules $V \subset W$
 that are isomorphic to $A$ (``type'')
 and admit
a complement $U \subset V$ isomorphic to $B$ (``cotype''). 
Then  
$$
(\XX^+,\XX^-)=(\Gras_A^B(W),\Gras_B^A(W)),
$$ 
with transversality of $(U,V) \in \XX^+ \times \XX^-$ 
being usual complementarity of 
subspaces ($W=U \oplus V$), is an affine pair geometry over $\K$:
it is a standard exercise in linear algebra to show
 that the set $U^\top$ of complements of $U$ carries
a natural structure of affine space.
Note that by our definition of Grassmannians as {\em complemented}
Grassmannians this affine space is never
empty; if $\K$ is a field this condition is automatic. 
(However, if $\K$ is a {\em topological} field or ring and
$W$ a topological $\K$-module, then one will prefer modified definitions
of ``restricted Grassmannians''
 by imposing conditions of closedness, of boundedness,
of Fredholm type or other, cf., e.g., \cite{PS86}, Chapter 7.
All such conditions lead to subgeometries of the algebraically defined
geometries considered here.)

\ssk
Fixing the decomposition $W=A \oplus B$ as base point $(o^+,o^-)$
in $(\XX^+,\XX^-;\top)$, by elementary linear algebra one may identify
the 
pair of linear spaces
 $((o^+)^\top,(o^-)^\top)$ with the pair 
$(V^+,V^-)=(\Hom(B,A),\Hom(A,B))$.
By still elementary (though less standard) linear algebra, one can now
give explicit expressions for the structure maps introduced in Section 2.2 
and check their fundamental identities, thus describing a direct link
between the geometry and its associated Jordan pair (see \cite{Be04}):
namely, elements of $\XX^+$ are realized as images of injective maps
$f:A \to W$, modulo equivalence under the general linear group
$\Gl_\K(A)$ ($f \sim f'$ iff $\exists g \in \Gl_\K(A)$:
$f'=f \circ g$), and similarly elements of $\XX^-$ are realized as kernels
of surjective maps $\phi:W \to A$, again modulo equivalence under
$\Gl_\K(A)$. Two such elements are transversal, $[f] \top [\phi]$, if and
only if $\phi \circ f:A \to A$ is a bijection, and the structure map 
$\PPP_r$ is now given by
$$
 \PPP_r([f],[\phi],[h])
= \big[ (1-r) f \circ (\phi \circ f)^{-1} + r h \circ (\phi \circ h)^{-1}
\big] .
$$
As in ordinary projective geometry, an affinization is given by writing
$f:A \to W = B \oplus A$ as ``column vector'' and normalizing the
second component to be ${\bf 1}_A$, the identity map of $A$, and similarly
for the ``row vector'' $\phi:B \oplus A \to A$ (see \cite{Be04}).
To get a feeling for the kind of non-linear formulas that appear in such
contexts, the reader may re-write the preceding formula for $\PPP_r$
 by replacing
$f,\phi$ and $h$ by such column-, resp.\ row-vectors, and then 
re-normalize the right-hand side, in order to get the formula for the
multiplication map in the affine picture. The special case $r=\frac{1}{2}$
(``midpoint map'') is particularly important from a Jordan-theoretic
point of view.

\ssk
One may as well also consider 
the {\em total} Grassmannian geometry of all complemented subspaces of $W$,
  $$
 \Gras(W):= \{ V \subset W | \, (V \, \mbox{submodule}), \,
\exists V' : W = V \oplus V'  (V' \, \mbox{submodule}) \} .
$$
Then the pair
$(\Gras(W),\Gras(W))$ still is a generalized projective geometry,
but it is not {\em connected} in general
(there is a natural notion of {\em connectedness} for any affine pair
geometry, see \cite{Be02}). 
In infinite dimension, the geometries $(\XX^+,\XX^-)$ introduced
above need not be connected neither, but in finite dimension over a field
they are; in fact, they then reduce to the usual Grassmannians
$\Gras_k(\K^n) \cong \Gl(n,\K)/(\Gl(k,\K) \times \Gl(n-k,\K))$
of $k$-spaces in $\K^n$. 

\ssk
The case $n=2k$, or, more generally, $A \cong B$, deserves special 
attention: in this case,
$\XX^+$ and $\XX^-$ really are the same sets; the identity map
$\XX^+ \to \XX^-$ is a canonical  antiautomorphism of the geometry.
Its properties are similar to the well-known {\em null-systems}
from classical projective geometry. In particular, in the case
of the projective line $(\Gras_1(\K^2),\Gras_1(\K^2))$,
this really is the canonical null-system coming from the (up to a scalar)
unique symplectic form on $\K^2$.

\msk
\nin {\bf (2)
Loop geomeries belong to loop algebraic structures.}
We fix some generalized projective geometry $\XX=(\XX^+,\XX^-)$, and
let $M$ be any set. Then the loop geometry
$$
(\Fon(M,\XX^+),\Fon(M,\XX^-))
$$
with ``pointwise transversality and structure maps'',
 is again a generalized projective geometry.
By arguments similar as above, we see that it corresponds to the
loop Jordan pair introduced in Example (2) of Section 1.5.
As explained there, the case of usual functions corresponds
to $V^+ = V^- = \Fon(M,\K)$ with pointwise triple product 
$(T(f,g,h))(p) = 2 f(p) \, g(p) \, h(p)$.
On the geometric side, then
$\XX^+ = \XX^- = \K \PP^1$ is the projective line, and
$\Fun(M,\XX^+)=\Fun(M,\XX^-)$ is the space of functions from $M$
into the projective line.

\msk \nin {\bf (3) The classical examples.} We briefly characterize the
geometries $(\XX^+,\XX^-)$ corresponding to the examples of item (3) of
section 1.5:

\begin{description}
\item{$\bullet$} 
full asociative algebras $A$ correspond to the {\em projective line
$(A \PP^1,A \PP^1)$ over $A$} (cf.\ \cite{BeNe04});
\item{$\bullet$} Hermitian elements of an involutive associative algebra
$(A,*)$ correspond to the $*$-Hermitian projective line over $A$
(cf.\ \cite{BeNe04}, Section 8.3), and similarly for
skew-Hermitian elements. {\em Lagrangian geometries}
 are a special case of this construction, corresponding to (skew-) Hermitian
or (skew-) symmetric matrices or operators 
\item{$\bullet$} conformal geometries (or ``spin factors''): here,
$\XX^+ \cong \XX^-$ is the projective quadric 
of $\K \oplus V \oplus \K$ with the same quadratic form
as in Section 1.5.
\end{description}

\nin 
As to the exceptional Jordan systems, their geometries are among those
constructed by Tits; to our knowledge, their structure of generalized
projective geometry has not yet been fully investigated.

\section{The universal model}

We have given above the construction of the generalized projective geometry
$(\XX^+,\XX^-)$ associated to a $3$-graded Lie algebra $\g$ via
homogeneous spaces $(G/P^-,G/P^+)$.
For a more detailed study, a geometric realization of these spaces is useful,
which in some sense generalizes the example of the Grassmannian geometries.
This ``universal model'' has been introduced in \cite{BeNe04} and  used
in \cite{BeNe05} to define (under suitable topological assumptions)
a manifold structure on $(\XX^+,\XX^-)$.

\subsection{Ordinary flag geometries}

We fix some $\K$-module $E$ and denote by $\FF_k$ the set of all flags
$$
\f = \big(0=\f_0 \subset \f_1 \subset \ldots \subset \f_k = E \big)
$$
of length $k$ in $E$ (i.e., the $\f_i$ are linear subspaces of $E$, and all 
inclusions are supposed to be strict).
We say that two flags $\e,\f \in \FF$ are {\em transversal} if they are
``crosswise complementary'' in the sense that
$$
\forall i = 1,\ldots,k: \quad \quad E = \f_i \oplus \e_{k-i}.
$$
It is a nice exercise in
linear algebra to show that $\e$ and $\f$ are transversal if,
and only if, they come from a {\it grading of $E$ 
of length $k$}:
given a grading $g=(\g_1,\ldots,\g_k)$, i.e., 
a decomposition
$
E = \g_1 \oplus \ldots \oplus \g_k,
$
we define two flags  $\f^+(g)$ and $\f^-(g)$ via
$$
 \f_i^+(g): = \g_1 \oplus \ldots \oplus \g_i, \quad \quad \f_{k-i}^-(g) : =
\g_{i+1} \oplus \ldots \oplus \g_k .
$$
Then $\f^+(g)$ and $\f^-(g)$ are obviously transversal,
but the converse is also true: every pair of transversal flags can be obtained
in this way!
From this exercise, one deduces easily that 
we get a linear pair geometry
$(\XX^+,\XX^-; \top)$ by taking 
$$
\XX^+ := \XX^- := \XX := \{ \f \in \FF_k | \, \, \exists \e \in \FF_k : \e \top \f \},
$$
the space of all flags that admit at least one transversal flag
(cf.\ \cite{BL06}).
For $k=2$, this is the total Grassmannan geometry, and similarly 
as in that case, one
may for general $k$
 single out connected components by taking {\em flags of fixed
type and cotype}. 
In general, these geometries will not be {\em affine} pair geometries, and
we do not know what kind of ``laws'' (in a sense generalizing the
laws of generalized projective geometries) can be used
to describe them. 

\subsection{Filtrations and gradings of Lie algebras}

Let us assume now that $\g:=E$ is a {\it Lie algebra}
with Lie bracket denoted by $[\cdot,\cdot]$,
and that all our gradings and filtrations are compatible with the Lie bracket in
the sense that
$$
[\f_i,\f_j] \subset \f_{i+j}, \quad [\g_i,\g_j] \subset \g_{i+j}.
$$
We will assume that  our Lie algebra is
{\em $2k+1$-graded}, i.e., 
as index set we take $I_k = \{ -k,-k+1, \ldots ,k \}$, and we
 say that a $2k+1$-grading is {\em inner} if it can be defined by an
{\em Euler-operator}, i.e., by an element $E \in \g$ such that
$[E,X] = i X$ for all $X \in \g_i$, $i \in I_k$.
We denote by
$\GG_k$ the set of all inner $2k+1$-gradings of $\g$ 
(this set can be identified
with set of all Euler operators).
An {\it inner $2k+1$-filtration of $\g$} is a flag
$$
\f = (0 = \f_{k+1} \subset \f_k \subset \ldots \subset \f_{-k}=\g)
$$
such that there is some inner $2k+1$-grading with
$\f_i = \g_i \oplus \g_{i+1} \oplus \ldots \oplus \g_k$, for all $i \in I_k$;
we write $\FF_k$ for the set of all such filtrations.
Then the following holds (\cite{BeNe04} for the case of $3$-gradings ($k=1$),
and \cite{Ch07} for the general case).
We have to impose some mild restrictions on the characteristic of $\K$ since we
want to use, for $X \in \f_1$, the operator $\exp(X)= e^{\ad(X)} =
\sum_{j=0}^\infty \frac{1}{j!} \ad(X)^j$, the sum being finite since
$[\f_1,\f_i] \subset \f_{i+1}$, so that $\ad(X)$ is nilpotent.

\bigskip \nin {\bf Theorem.}  {\em
Assume that $\K$ is a commutative base ring such that the integers
$2,3,\ldots, 2k+1$ are invertible in $\K$. 
\begin{description}
\item{\rm (1)} Two inner $2k+1$-filtrations $\f$ and $\e$ are transversal if, and only
if, there exists a $2k+1$-grading $g = (\g_{-k},  \ldots, \g_k)$ such
that
$\f = \f^+(g)$ and $\e = \f^-(g)$.
\item{\rm (2)} The geometry of inner $2k+1$-filtrations $(\FF_k,\FF_k;\top)$
is a linear pair geometry. More precisely, for any $\f \in \FF_k$, the
nilpotent Lie algebra $\f_1$ is in bijection with $\f^\top$ via
$X \mapsto \exp(X).\e$, for an arbitrary choice of $\e \in \f^\top$.
\end{description}
}

\nin
Now assume that $k=1$, i.e., we are in the $3$-graded case.
Then $[\f_1,\f_1] \subset \f_2 = 0$ (i.e., $\f_1$
is an abelian subalgebra), whence Part (2) of the
theorem says that we have a simply transitive action of the vector group
$\f_1 \cong \exp(\f_1)$ on $\f^\top$, and thus we see that the geometry of
inner $3$-filtrations is an {\em affine} pair geometry.
In fact,  fixing some inner $3$-grading as base point in the space of all
inner $3$-gradings, it now easily follows that the geometry $(G/P^-,G/P^+)$
constructed in the preceding chapter is imbedded into the geometry from the 
theorem simply by taking the $G$-orbit of the base point.

\ssk
As mentioned above, this ``universal
geometric realization'' turns out to be useful for giving a precise description
of the intersections of the ``chart domains''
 $\alpha^\top$ and  $\beta^\top$ for $\alpha,\beta \in \XX^-$
 and for calculating the
corresponding ``transition functions'' (see \cite{BeNe04}); 
 calculations become similar to the case of Grassmannian geometries (Example
(1) in Section 2.3 above), the difference being that the usual fractional
linear transformations have to be replaced by certain {\em fractional
quadratic} transformations (this difference
corresponds to the fact that we are now
working with true flags instead of single subspaces). 
Using these purely algebraic results, one can now
give necessary and sufficient conditions for defining on 
$(\XX^+,\XX^-)$ the structure
of a smooth manifold (\cite{BeNe05}). 
The theory works nicely over general topological base fields and rings
(cf.\ \cite{Be06} for basic differential geometry and Lie theory in this
general framework); complex or real Banach, Fr\'echet or locally
convex manifolds are special cases of it.
For higher gradings, similar 
results can be expected (J.\ Chenal, work in progress, \cite{Ch07}).

\section{The geometry of states}

Let us come back to the analogies, mentioned in the introduction,
 of the preceding constructions with methods of commutative and non-commutative
geometry. In the language of classical physics, one wants to recover the
pure states of a system (the point space $M$) from its observables (the
function algebra $\cA =\Reg(M,\R)$). The usual procedure
 is to look at the (maximal)
  ideals $I_p = \{ f \in \cA | \; f(p)=0 \}$ corresponding to points 
$p \in M$. Left or right ideals in associative algebras are generalized
by {\em inner ideals} in Jordan theory. Therefore we shall interprete
{\em states} in generalized projective geometries as the geometric objects
corresponding to inner ideals.

\subsection{Intrinsic subspaces}

Assume $(\XX^+,\XX^-;\top)$ is a linear pair geometry over the commutative ring
$\K$. 

\begin{definition}
A pair $(\YY^+,\YY^-)$ of subsets
 $\YY^\pm \subset \XX^\pm$  is called a
{\em subspace of $(\XX^+,\XX^-)$}, if
for all $x \in \YY^\pm$ there exists an element  $\alpha \in \YY^\mp$ 
such that $x \top \alpha$, and if
for every pair $(x,\alpha) \in (\YY^+ \times \YY^-)^\top$ the set
$\YY_\alpha:=\YY \cap \alpha^\top$ is a {\it linear} 
subspace of $(\alpha^\top,x)$ and
$\YY_x':=\YY' \cap x^\top$ is a {\it linear} subspace of $(x^\top,\alpha)$.

\ssk
A subset
 $\II \subset \XX^+$ is called a
{\em state} or
 {\em intrinsic subspace (in $\XX^+$)}, 
if  $\II$ ``appears linearly to {\it all possible} observers'':
for all  $\alpha \in \XX^-$ with
$\alpha^\top \cap \II \not= \emptyset$ and for all 
$x \in \alpha^\top \cap \II$, the
set $(\alpha^\top \cap \II,x)$ is a linear subspace of $(\alpha^\top,x)$.
\end{definition}

For instance, if $(\XX^+,\XX^-)=(\Gras_1(\K^{n+1}),\Gras_n(\K^{n+1}))$
is an ordinary projective geometry, then all projective subspaces of 
$\XX^+$ are intrinsic subspaces, and every affine subspace of an
affinization of $\XX^+$ is obtained in this way.
In contrast, for Grassmannian geometries of higher
rank the situation becomes more complicated: only rather specific linear 
subspaces of the affinization $M(p,q;\K)$ of $\Gras_p(\K^{p+q})$ are obtained
in this way, namely the so-called {\em inner ideals}. This observation
generalizes to all geometries
 $(\XX^+,\XX^-)$ associated to a Jordan pairs $(V^+,V^-)$:
subspaces containing the base point $(o^+,o^-)$ correspond to subpairs
of $(V^+,V^-)$, and intrinsic subspaces of $\XX^+$ containing $o^+$ to
{\em inner ideals} $I \subset V^+$, i.e., submodules such that
$T^+(I,V^-,I) \subset I$ (see \cite{BL06}).
The notions of {\em minimal, maximal, principal,...} inner ideals may also
be suitably translated into a geometric language (minimal intrinsic subspaces
will also be called {\em intrinsic lines} or {\em pure states}).
The collection of  states associated to a linear pair geometry can 
again be turned into a pair geometry by introducing the following
transversality relation:

\begin{definition}
Let $\II$ be an intrinsic subspace in $\XX^+$ and $\JJ$ one in $\XX^-$.
We say that $\II$ and $\JJ$ are {\em transversal}, and we write again
$\II \top \JJ$, if

\begin{description}
\item{(1)} the pair $(\II,\JJ)$ is a subspace,
\item{(2)} the linear pair geometry $(\II,\JJ)$ is {\em faithful} in the following
 sense:
\end{description}

A linear pair geometry $(\YY^+,\YY^-)$ is called {\em faithful} if 
$\YY^-$ is faithfully represented by its effect of
 linearizing $\YY^+$, and vice versa: whenever 
 $\alpha^\top = \beta^\top$ as sets {\em and} as linear spaces
(with respect to some origin $o$), then $\alpha=\beta$, and the dual
property holds.  
\end{definition}

\nin
Let us give some motivation for this definition (which does not appear in
\cite{BL06}). First of all, in general, a linear pair geometry need not
be faithful (take, for instance, a pair of $\K$-modules with trivial
structure maps), but ordinary projective geometries over a field are.
(In case of geometries corresponding to Jordan pairs,
faithfulness corresponds to {\em non-degeneracy} in the Jordan-theoretic
sense.) Next note that,
even if the geometry $(\XX^+,\XX^-)$ was faithful, the geometry
$(\II,\XX^-)$ will in general not be faithful: there is a {\em kernel},
i.e., the equivalence relation ``$\alpha \sim_\II \beta$ iff
$\alpha$ and $\beta$ induce the {\em same} linear structure on $\II$''
will be non-trivial on $\XX^-$.
The condition $\II \top \JJ$ means, then, that the space $\JJ \subset \XX^-$
is transversal to the fibers of this equivalence relation, in the usual
sense. 
(This is the geometric translation of the notion of {\em complementation of
inner ideals} introduced in
\cite{LoNe95}: the fiber of the equivalence relation $\sim_\II$ corresponds
to  the {\em kernel of an inner ideal}  introduced in \cite{LoNe95}).
For instance, if $\II$ is a projective line in an ordinary
 projective space $\XX^+ =\Gras_1(\K^{n+1})$, corresponding to a
$2$-dimensional subspace $I$ of $\K^{n+1}$, then $\JJ$ should be a
projective subspace of the dual projective space $\XX^-=\Gras_n(\K^{n+1})$
such that different elements of $\JJ$ define different affine lines
in $\II$. If $\JJ$ is the set of hyperplanes lying over some $k$-dimensional
subspace of $\K^{n+1}$, then this means that $I$ should contain some
complement of $J$. Vice versa, $J$ then also should contain some complement
of $I$, and hence $I$ and $J$ have to be complements of each other in
$\K^{n+1}$. Thus $(\II,\JJ)$ is a projective line which is faithfully
imbedded in the projective geometry. Moreover, in this example we
clearly recover the Grassmannian geometry of $2$-spaces as the
geometry of intrinsic lines in the $1$-Grassmannian.

\ssk Parts of the observations from the example generalize:
let $\SSSS^\pm$ be the set of all intrinsic subspaces
 $\II$ in $\XX^\pm$ that admit 
a transversal intrinsic subspace $\JJ$ in $\XX^\mp$.
It is a natural question to ask whether $(\SSSS^+,\SSSS^-)$ is again a
linear pair geometry, or under which conditions on $(\XX^+,\XX^-)$ and
on the states we obtain one. 
For the time being, we have no general answers,
but results from \cite{BL06} point into the direction that, indeed, 
for the most interesting cases we get geometries $(\SSSS^+,\SSSS^-)$ that are
again linear pair geometries. More precisely, if $(\XX^+,\XX^-)$ corresponds to
a Jordan pair $(V^+,V^-)$ and we consider only states
that are associated, via a Peirce-decomposition, to idempotents, we
obtain geometries coming from $5$-graded Lie algebras  (\cite{BL06},
Theorem 5.8), and these
 are indeed linear pair geometries, according to Theorem 3.2.


\subsection{Examples}

{\bf (1)} Classical states and observables revisited.
Let $M$ be a set and $\cA = \Fun(M,\K)$ the Jordan algebra of all fonctions
from $M$ to $\K$. It corresponds to the geometry $(\XX,\XX)$ with
$\XX= \Fun(M,\K \PP^1)$. 
We claim that (if $\K$ is a field)
 the set $M$ can be naturally identified with the
geometry of all pure states (intrinsic lines)
running  through the ``zero function'' $f_0 \equiv o$,
 by associating to a point $p \in M$ the intrinsic line
$$
L_p = \{ f:M \to \K \PP^1 | \, f(x) = o \, \, \mbox{if} \, \, x \not= p \}
\subset \Fon(M,\K \PP^1) .
$$
Indeed, note first that $L_p$ is indeed a minimal intrinsic subspace
because the geometry $(\K \PP^1,\K \PP^1)$ itself does not admit
any proper intrinsic subspaces that are not points.
Now
assume $f:M \to \K \PP^1$ is a function taking non-zero values at two
different points $p,q \in M$.
Then a look at the direct product geometry 
$(\K \PP^1,\K \PP^1) \times (\K \PP^1,\K \PP^1)$ shows that
the intrinsic subspace generated by $f$ contains at least all functions
obtained from $f$ by altering the values at $p$ and $q$ in an arbitrary way,
hence contains $\K \PP^1 \times \K \PP^1$ and thus is not  minimal,
proving our claim.
If $M$ is a smooth manifold and we are more interested in loop
geometries of smooth functions, we may also work with  {\it maximal}
intrinsic subspaces, which then are spaces of smooth
 functions vanishing at a point.
We thus mimick constructions 
that imbed ``classical geometry''
into commutative algebra. 
%

\msk \nin
{\bf (2)} 
Grassmann geometries $\XX^+ = \Gras_A^B(W)$.
Let us fix a flag in $W$
 of length two, $\f: 0 \subset \f_1 \subset \f_2 \subset W$.
Then the set of all subspaces of $W$ that are ``squeezed'' by this
flag, 
$$
\II_\f := \{ E \in \XX^+ | \, \f_1 \subset E \subset \f_2 \},
$$
is an intrinsic subspace of $\XX$, and if $\K$ is a field and
$W$  is finite-dimensional over $\K$,
then all intrinsc subspaces are of this form for a suitable flag $\f$
(see \cite{BL06}).
Moreover, one can show that
the intrinsic subspaces $\II_\f$ and $\II_\e$ are transversal in
the sense defined above if, and only if, 
the flags $\e$ and $\f$ are transversal
in the sense defined in Section 3.1 above.
Therefore the geometry of intrinsic subspaces of the Grassmannian geometry
is a geometry of flags of length two and hence is a linear pair geometry.
The case of ordinary projective geometry (i.e., $A \cong \K$)
is somewhat degenerate: in this case there is not much choice for the
first component $\f_1$ of the flag $\f$, and therefore the intrinsic
subspace $\II_\f$ is already determined by $\f_2$ alone; 
we get back the usual projective subspaces of a projective space.
In all other cases, the geometry of intrinsic subspaces 
is a true flag geometry corresponding
to $5$-gradings  and hence is a linear, but no longer an affine pair geometry.

\msk \nin
{\bf (3)} 
Quantum states.
In the language of quantum mechanics (widely used also
in non-commutative geometry), {\em states} are
defined as positive normalized linear functionals $\phi: \cA \to \C$
on a $C^*$-algebra $(\cA,*)$.
They form a convex set, whose extremal points are interpreted as
{\em pure states}. These, in turn, admit a Jordan theoretic interpretation
via {\it (primitive) idempotents} -- see \cite{FK94} for the theory
of finite-dimensional {\em symmetric cones} and their geometry; there is
a vast litterature on infinite dimensional generalizations, cf., e.g.,
\cite{HOS84}, \cite{ER92}, \cite{Up85}. 
Whereas the notion of state
just mentioned depends highly on the ordered structure
 of the base field $\R$ (via {\it positivity} and {\em convexity}),
this is not the case for the notion of inner ideal and intrinsic subspace,
which hence seem to be more general and more geometric.
Accordingly, it has been advocated to use inner ideals as a basic ingredient
for an approach to quantum mechanics (see \cite{F80}).

%

%

\vfill\eject

\end{document}